\documentclass{amsart}

\usepackage{amssymb}
\usepackage{amsthm}
\usepackage{amsmath}
\usepackage{graphicx}
\usepackage{enumitem}

\theoremstyle{plain}
\newtheorem*{thm*}{Theorem}
\newtheorem{thm}{Theorem}
\newtheorem{prop}[thm]{Proposition}
\newtheorem{lem}[thm]{Lemma}
\newtheorem{cor}[thm]{Corollary}

\theoremstyle{definition}
\newtheorem{defn}[thm]{Definition}
\newtheorem{claim}{Claim}[thm]

\newenvironment{statement}[1]{\begin{quote}\textup{{#1}:} \,}{\end{quote}}

\newcommand{\id}{\mathrm{id}}
\newcommand{\diam}{\mathrm{diam}}
\newcommand{\proj}{\mathrm{proj}}
\newcommand{\length}{\mathsf{len}}
\newcommand{\interior}{\mathrm{int}}
\newcommand{\dsup}{d_{\mathrm{sup}}}
\newcommand{\lam}{\mathcal{L}}

\renewcommand{\mod}{\mathrm{mod}}

\begin{document}

\title[A canonical parameterization of paths in $\mathbb{R}^n$]{A canonical parameterization\\of paths in $\mathbb{R}^n$}
\author{L. C. Hoehn \and L. G. Oversteegen \and E. D. Tymchatyn}
\date{September 9, 2016}

\address[L.\ C.\ Hoehn]{Nipissing University, Department of Computer Science \& Mathematics, 100 College Drive, Box 5002, North Bay, Ontario, Canada, P1B 8L7}
\email{loganh@nipissingu.ca}

\address[L.\ G.\ Oversteegen]{University of Alabama at Birmingham, Department of Mathematics, Birmingham, AL 35294, USA}
\email{overstee@uab.edu}

\address[E.\ D.\ Tymchatyn]{University of Saskatchewan, Department of Mathematics and Statistics, 106 Wiggins road, Saskatoon, Canada, S7N 5E6}
\email{tymchat@math.usask.ca}

\thanks{The first named author was partially supported by NSERC grant RGPIN 435518, and by the Mary Ellen Rudin Young Researcher Award}
\thanks{The second named author was partially supported by NSF grant DMS-0906316}
\thanks{The third named author was partially supported by NSERC grant OGP-0005616}

\subjclass[2010]{Primary 54F15; Secondary 54C20, 54F50}
\keywords{path, length, parameterization of arcs, lamination}

\begin{abstract}
For sufficiently tame paths in $\mathbb R^n$, Euclidean length provides a canonical parametrization of a path by length.  In this paper we provide such a parametrization for all continuous paths.  This parametrization is based on an alternative notion of path length, which we call $\length$.  Like Euclidean path length, $\length$ is invariant under isometries of $\mathbb R^n$, is monotone with respect to sub-paths, and for any two points in $\mathbb R^n$ the straight line segment between them has minimal $\length$ length.

Unlike Euclidean path length, the $\length$ length of \emph{any} path is defined (i.e., finite) and $\length$ is continuous relative to the uniform distance between paths.  We use this notion to obtain characterizations of those families of paths which can be reparameterized to be equicontinuous or compact.  Finally, we use this parametrization to obtain a canonical homeomorphism between certain families of arcs.
\end{abstract}

\maketitle

\section{Introduction}
\label{sec:introduction}

A \emph{path} in $\mathbb{R}^n$ is a continuous function $\gamma$ from a closed interval $[a,b] \subset \mathbb{R}$ to $\mathbb{R}^n$.  Given $z_1,z_2 \in \mathbb{R}^n$, denote by $\overline{z_1 z_2}$ the straight line segment path $t \mapsto (1-t) z_1 + t z_2$, $t \in [0,1]$.

Given a path $\gamma: [a,b] \to \mathbb{R}^n$, the \emph{Euclidean path length} of $\gamma$, denoted $L_E(\gamma)$, is defined by the formula
\[ L_E(\gamma) = \sup \left\{ \sum_{i=1}^n |\gamma(x_{i-1}) - \gamma(x_i))|: a = x_0 < x_1 < \cdots < x_n = b, \; n \in \mathbb{Z}^+ \right\} \in [0,\infty] ,\]
where $|z_1 - z_2|$ denotes the Euclidean distance between points $z_1,z_2 \in \mathbb{R}^n$.

If a sequence of smooth paths $\gamma_i$ converges to $\gamma_\infty$ in $C^\infty$ (in the sense that the paths $\gamma_i$ and their derivatives $\gamma_i'$ converge uniformly to $\gamma_\infty$ and $\gamma_\infty'$, respectively), then $L_E(\gamma_i) \to L_E(\gamma_\infty)$ and, hence, path length provides a canonical parameterization of this entire family.  One of the main goals in this paper is to extend such results to the topological category.  For this reason we introduce a new notion of path length which is defined for all paths and behaves well with respect to uniform convergence of paths.

\bigskip
The function $L = L_E$ satisfies the following basic properties for the path $\gamma: [a,b] \to \mathbb{R}^n$:
\begin{enumerate}[label=(\textbf{L\arabic{*}})]
\item \label{enum:L1} If $A \subset [a,b]$ is a closed subinterval, then $L(\gamma {\upharpoonright}_A) \leq L(\gamma)$;
\item \label{enum:L2} If $c \in (a,b)$, then $L(\gamma) = L(\gamma {\upharpoonright}_{[a,c]}) + L(\gamma {\upharpoonright}_{[c,b]})$;
\item \label{enum:L3} If $\Phi: \mathbb{R}^n \to \mathbb{R}^n$ is an isometry, then $L(\Phi \circ \gamma) = L(\gamma)$;
\item \label{enum:L4} $L(\gamma) =$ the supremum, taken over all partitions $a = x_0 < x_1 < \cdots < x_n = b$ of $[a,b]$, of the values $L(P)$ where $P$ is the polygonal path with vertices $\gamma(x_0),\ldots,\gamma(x_n)$;
\end{enumerate}
and moreover, we have
\begin{enumerate}[label=(\textbf{L\arabic{*}})]
\setcounter{enumi}{4}
\item \label{enum:L5} $L(\overline{0 e}) = 1$, where $e = (1,0,\ldots,0) \in \mathbb{R}^n$.
\end{enumerate}

Conversely, any function $L$ defined on the set of all paths which satisfies the properties \ref{enum:L1} through \ref{enum:L5} must be equal to $L_E$ (and any function $L$ which satisfies the properties \ref{enum:L1} through \ref{enum:L4} must be a scalar multiple $c \cdot L_E$ of $L_E$, where $c = L(\overline{0 e})$).  Indeed, one can use properties \ref{enum:L1}, \ref{enum:L2}, \ref{enum:L3}, and \ref{enum:L5} to show that $L(\overline{a b}) = b - a$ for any $a,b \in \mathbb{R} \subset \mathbb{R}^n$ with $a < b$.  Then by properties \ref{enum:L2} and \ref{enum:L3} it follows that the length of any polygonal path is equal to the sum of the (Euclidean) distances between consecutive vertices.  We then conclude by \ref{enum:L4} that $L(\gamma) = L_E(\gamma)$ for all paths $\gamma$.

There are a number of results in metric geometry pertaining to when a given metric on a Euclidean space is equal to the Euclidean metric; \cite{BlumenthalMenger1970}, \cite{Busemann1955}, and \cite{Menger1928} each survey a variety of such results.  Much of this work is related to Hilbert's fourth problem.  The length function introduced in this paper contributes to the corresponding program for path length functions by illustrating that there are other length functions which have many properties in common with the Euclidean length.

In light of the above discussion, to provide a genuinely different path length function from the Euclidean length, one must give up at least one of the properties \ref{enum:L1} through \ref{enum:L4}.  In Section \ref{sec:definition}, we define a path length function, called ``$\length$'', such that $L = \length$ satisfies properties \ref{enum:L1}, \ref{enum:L3}, and \ref{enum:L4} (see Propositions \ref{prop:basic properties 1}\ref{enum:iso inv}, \ref{prop:basic properties 2}\ref{enum:monotone}, and \ref{prop:length continuous} below), as well as the following weaker form of \ref{enum:L2} (see Proposition \ref{prop:basic properties 2}\ref{enum:subadditive}):
\begin{enumerate}[label=(\textbf{L\arabic{*}$'$})]
\setcounter{enumi}{1}
\item \label{enum:L2'} If $c \in (a,b)$, then $L(\gamma) \leq L(\gamma {\upharpoonright}_{[a,c]}) + L(\gamma {\upharpoonright}_{[c,b]})$;
\end{enumerate}

Furthermore, this length function has the following additional properties not enjoyed by the Euclidean length $L_E$:
\begin{itemize}
\item $\length(\gamma) < 1$ for any path $\gamma$;
\item $\length$ is continuous as a function from the space of all maps $[0,1] \to \mathbb{R}^n$ (with the uniform metric) to $\mathbb{R}$.
\item $\length$ is defined for any continuous function $\gamma$ from a locally connected continuum $X$ to $\mathbb{R}^n$;
\end{itemize}
Moreover, this length function can differentiate between paths whose Euclidean lengths are infinite.  For instance, if $\gamma: [a,b] \to \mathbb{R}^n$ is a path and $[c,d]$ is a subinterval of $[a,b]$ such that $\gamma$ is non-constant on at least one component of $[a,b] \smallsetminus [c,d]$, then $\length(\gamma) > \length(\gamma {\upharpoonright}_{[c,d]})$, even if both of these paths have infinite Euclidean length.

\bigskip
A very similar function is developed by Cannon et.\ al.\ in \cite{CannonConnerZastrow2002}, which is called the \emph{total oscillation} of a path.  The most notable difference is that the total oscillation is not invariant under isometries of $\mathbb{R}^n$.

Another similar function is given by Morse in \cite{Morse1936}, called the \emph{$\mu$-length}, which is defined for paths into any metric space.

\bigskip
After establishing the above properties in Section \ref{sec:properties}, we use $\length$ in Section \ref{sec:parameterization length} to obtain a standard parameterization of all paths in $\mathbb{R}^n$.  This yields characterizations of those families of paths which may be reparameterized so as to be equicontinuous or compact.  These results extends classical results on families of paths having finite Euclidean length.  In Section \ref{sec:parameterization midpoint}, we develop a second canonical parameterization of all paths in $\mathbb{R}^n$ with all of the above properties, and which also commutes with a reversal of orientation of a path.  This second parameterization yields a canonical extension of a bijection between the endpoint sets of two arcs to a homeomorphism between the two arcs.  We use this notion to construct homeomorphisms between certain families of pairwise disjoint arcs in Theorem \ref{thm:lamination}.

The results of this paper have already been used in two other papers.  In \cite{OT10} it was shown that any isotopy of a planar continuum can be extended to an isotopy of the entire plane.  Using Theorems \ref{thm:lamination} and \ref{thm:lamination cont}, this result is extended in \cite{HOT14} to a more general class of planar compacta.  It was shown in \cite{BR89} that any two points in a closed topological disk $D$ in the plane can be connected by a unique arc $A$ in $D$ which has the property that any subarc of $A$ which connects two points, neither of which is an endpoint of $A$, has minimal (finite) Euclidean length among all such arcs.  In \cite{HOT14a} this result is generalized to shortest paths (in the sense of $\length$ length, and in the above Euclidean sense for proper subpaths) in the closure of any homotopy class in an open connected subset of the plane with arbitrary boundary.


\section{Definition of the function $\length$}
\label{sec:definition}

A \emph{generalized path} is a continuous function $\gamma: X \to \mathbb{R}^n$, where $X$ is a locally connected metric continuum.

Given $n \geq 2$, there is a length function $\length_n$ defined for generalized paths $X \to \mathbb{R}^n$.  In this paper, to simplify the definition and arguments below, we restrict our attention to the case $n = 2$, and give a definition of $\length = \length_2$.  The case $n > 2$ proceeds similarly; the primary differences being that we cut by $(n-1)$-dimensional hyperplanes instead of lines (see below), and the parameter $t$ below varies through the $(n-1)$-dimensional real projective space instead of $[0,1]$ (which we use to parameterize the semi-circle $\{e^{t \pi i}: t \in [0,1]\}$).

The reader may find it easier to work through this construction for an ordinary path $\gamma: [a,b] \to \mathbb{C}$ instead of a generalized path, on a first reading.

For notational convenience, we will identify $\mathbb{R}^2$ with $\mathbb{C}$.

\bigskip
For $j \in \mathbb{Z}$, let $S_j$ denote the closed horizontal strip $\{a + ib: a \in \mathbb{R}, b \in [j,j+1]\}$ in the plane $\mathbb{C}$.  Given $x,t \in [0,1]$, $\mu \in (0,1]$, and $j \in \mathbb{Z}$, let $S^{x,t,\mu}_j = \mu e^{t \pi i}(S_j + ix)$.  If $A \subset \mathbb{C}$, define $\|A\|_t = \diam(\proj^\perp_t(A))$, where $\proj^\perp_t$ denotes the orthogonal projection of $\mathbb{C}$ onto the line $\{r e^{(t+\frac{1}{2}) \pi i}: r \in \mathbb{R}\}$ and $\diam$ denotes the diameter in the Euclidean metric.

Fix a generalized path $\gamma: X \to \mathbb{C}$.

The following lemma will be used in the definition of the function $\length$ below.

\begin{lem}
\label{lem:components large projections}
For any $(x,t,\mu) \in [0,1] \times [0,1] \times (0,1]$ and any $\varepsilon > 0$, there are only finitely many components $C$ of the sets $\gamma^{-1}(S^{x,t,\mu}_j)$ ($j \in \mathbb{Z}$) with $\|\gamma(C)\|_t \geq \varepsilon$.
\end{lem}

\begin{proof}
We may assume that $\varepsilon \leq \frac{\mu}{2}$.  Suppose for a contradiction that there are infinitely many distinct components $\{C_n\}_{n=0}^\infty$ of the sets $\gamma^{-1}(S^{x,t,\mu}_j)$ ($j
\in \mathbb{Z}$) with $\|\gamma(C)\|_t \geq \varepsilon$.

For each $n$ let $j(n) \in \mathbb{Z}$ be the integer for which $\gamma(C_n) \subset S^{x,t,\mu}_{j(n)}$, and let $p_n \in C_n$ be such that $d(\gamma(p_n), \partial S^{x,t,\mu}_{j(n)}) = \varepsilon$, where $d$ denotes the Euclidean metric in $\mathbb{R}^n$.  Observe that by local connectivity of $X$, for each $n$ we have $\gamma(\partial C_n) \subset \partial S^{x,t,\mu}_{j(n)}$.

Let $p \in X$ be an accumulation point of the set $\{p_n\}_{n=0}^\infty$, and let $U$ be an open neighborhood of $p$ which is small enough so that $\diam(\gamma(U)) < \varepsilon$.  Then for any $n$ such that $p_n \in U$, we have $\gamma(U) \cap \partial S^{x,t,\mu}_{j(n)} = \emptyset$, hence $U \cap \partial C_n = \emptyset$, and so $U \cap C_n$ is closed and open in $U$.  It follows that $U$ cannot be connected, which is a contradiction since $X$ is locally connected.
\end{proof}

Given $x,t \in [0,1]$ and $\mu \in (0,1]$, let ${\langle C^{x,t,\mu}_n \rangle}_{n=0}^\infty$ enumerate the collection of all components of the sets $\gamma^{-1}(S^{x,t,\mu}_j)$ ($j \in \mathbb{Z}$) which have non-degenerate image under the map $\proj^\perp_t$, ordered so that $\|\gamma(C^{x,t,\mu}_n)\|_t \geq \|\gamma(C^{x,t,\mu}_{n+1})\|_t$ for all $n$ (this is possible by Lemma \ref{lem:components large projections}).

Define
\[ L^{x,t,\mu}(\gamma) = \sum_{n=0}^\infty \frac{\|\gamma(C^{x,t,\mu}_n)\|_t}{2^n} \]
and define the \emph{length of $\gamma$} by
\[ \length(\gamma) = \int_0^1 \! \int_0^1 \! \int_0^1 \! L^{x,t,\mu}(\gamma) \, dx \, dt \, d\mu .\]

If $X \subset \mathbb{C}$ is a locally connected continuum, define $\length(X) = \length(\id_X)$.

Observe that if $\sigma$ is any injective function of the non-negative integers to themselves, then
\begin{equation}
\tag{$\ast$} \sum_{n=0}^\infty \frac{\|\gamma(C^{x,t,\mu}_{\sigma(n)})\|_t}{2^n} \leq L^{x,t,\mu}(\gamma) .
\end{equation}

It remains to show that the function $L^{x,t,\mu}(\gamma)$ is in fact integrable, so that the above definition of the function $\length$ makes sense.  This is accomplished in Lemma \ref{lem:lower semicontinuity} below.

\begin{lem}
\label{lem:component subset}
Let $C$ be a component of $\gamma^{-1}(S^{x,t,\mu}_j)$ for some $x,t,\mu,j$ which has non-degenerate image under the map $\proj^\perp_t$, and let $\varepsilon > 0$.  Then there exists a subcontinuum $D \subset C$ such that $\gamma(D) \subset \interior(S^{x,t,\mu}_j)$ and $\|\gamma(D)\|_t \geq \|\gamma(C)\|_t - \varepsilon$.
\end{lem}

\begin{proof}
For the purposes of this argument, let us naturally identify $\mathbb{R}$ with the rotated line $\{r e^{(t+\frac{1}{2}) \pi i}: r \in \mathbb{R}\}$ which is the range of the map $\proj^\perp_t$.

Let $s_1,s_2 \in \mathbb{R}$ be such that $s_1 < s_2$ and $\proj^\perp_t(\gamma(C)) = [s_1,s_2]$ (and hence $\|\gamma(C)\|_t = s_2 - s_1$).  We may assume that $\varepsilon < \frac{s_2 - s_1}{2}$.  Let $S'$ denote the narrower (closed) strip $(\proj^\perp_t)^{-1}([s_1 + \frac{\varepsilon}{2}, s_2 - \frac{\varepsilon}{2}]) \subset \interior(S^{x,t,\mu}_j)$.  Then $C \cap \gamma^{-1}(S')$ must have a component $D$ such that $\proj^\perp_t(\gamma(D)) = [s_1 + \frac{\varepsilon}{2}, s_2 - \frac{\varepsilon}{2}]$ (see e.g.\ Theorem 5.2 of \cite {Nadler1992}).  This $D$ is as desired.
\end{proof}

A real-valued function $f$ is \emph{lower semicontinuous} if $f^{-1}((\alpha, \infty))$ is open for every $\alpha \in \mathbb{R}$.  Note that a lower semicontinuous function is Borel, hence (Lebesgue) integrable.

\begin{lem}
\label{lem:lower semicontinuity}
For a fixed generalized path $\gamma: X \to \mathbb{C}$, put $L(x,t,\mu) = L^{x,t,\mu}(\gamma)$.  Then the function $L(x,t,\mu)$ from $[0,1] \times[0,1] \times (0,1]$ to $\mathbb{R}$ is lower semicontinuous, hence integrable.
\end{lem}

\begin{proof}
Fix a number $\alpha \in \mathbb{R}$, and suppose $L^{x,t,\mu}(\gamma) > \alpha$.  Choose $N$ large enough so that $\sum_{n=0}^N \frac{\|\gamma(C^{x,t,\mu}_n)\|_t}{2^n} > \alpha$.

For each $n \in \{0,1,\ldots,N\}$ let $j(n)$ be such that $C^{x,t,\mu}_n$ is a component of $\gamma^{-1}(S^{x,t,\mu}_{j(n)})$.  Then, by Lemma \ref{lem:component subset}, for each $n$ we can find a proper subcontinuum $D_n \subset C^{x,t,\mu}_n$ such that $\gamma(D_n)$ is contained in the interior of $S^{x,t,\mu}_{j(n)}$, and so that
\[ \sum_{n=0}^N \frac{\|\gamma(D_n)\|_t}{2^n} > \alpha .\]

Let $\varepsilon_1 > 0$ be small enough so that if $|x' - x|, |t' - t|, |\mu' - \mu| < \varepsilon_1$, then $\gamma(D_n) \subset S^{x',t',\mu'}_{j(n)}$ for each $n \in \{0,1,\ldots,N\}$, and moreover
\begin{equation}
\label{eq:1}
\sum_{n=0}^N \frac{\|\gamma(D_n)\|_{t'}}{2^n} > \alpha.
\end{equation}

For each pair of numbers $n_1 < n_2$ in $\{0,1,\ldots,N\}$ with $j(n_1) = j(n_2)$, find an open set $A_{n_1,n_2} \subset X$ such that $C^{x,t,\mu}_{n_1} \subset A_{n_1,n_2} \subset \overline{A_{n_1,n_2}} \subset X \smallsetminus C^{x,t,\mu}_{n_2}$ and $\partial A_{n_1,n_2} \cap \gamma^{-1}(S^{x,t,\mu}_{j(n_1)}) = \emptyset$; that is, $\gamma(\partial A_{n_1,n_2}) \cap S^{x,t,\mu}_{j(n_1)} = \emptyset$.

Let $0 < \varepsilon_2 < \varepsilon_1$ be small enough so that if $|x' - x|, |t' - t|, |\mu' - \mu| < \varepsilon_2$, then $\gamma(\partial A_{n_1,n_2}) \cap S^{x',t',\mu'}_{j(n_1)} = \emptyset$ for every pair of numbers $n_1 < n_2$ in $\{0,1,\ldots,N\}$ with $j(n_1) = j(n_2)$.  Since $\partial A_{n_1,n_2}$ separates $D_{n_1}$ from $D_{n_2}$ in $X$, it follows that $D_{n_1}$ and $D_{n_2}$ are contained in distinct components of $\gamma^{-1}(S^{x',t',\mu'}_{j(n_1)})$.  Therefore, for such $x',t',\mu'$, by ($\ast$) and (\ref{eq:1}) we have
\[ L^{x',t',\mu'}(\gamma) \geq \sum_{n=0}^N \frac{\|\gamma(D_n)\|_{t'}}{2^n} > \alpha .\]

Thus, the set $\{(x,t,\mu): L^{x,t,\mu}(\gamma) > \alpha\}$ is open in $[0,1] \times [0,1] \times (0,1]$, and so $L(x,t,\mu)$ is a lower semicontinuous function.
\end{proof}

Thus the function $\length$ is well-defined.  Observe that the set $\gamma(C^{x,t,\mu}_n)$ is contained in some strip $S^{x,t,\mu}_j$ having width $\mu$, hence $\|\gamma(C^{x,t,\mu}_n)\|_t \leq \mu$. It follows that $L^{x,t,\mu}(\gamma) < 2\mu$, and therefore $\length(\gamma) < 1$.

\medskip

It can easily be seen that $\length(\overline{0x}) \to 1$ as $x \to \infty$, $x \in \mathbb{R}$.  It follows from Propositions \ref{prop:basic properties 1}\ref{enum:iso inv} and \ref{prop:basic properties 3} below that if $\gamma_m: X_m \to \mathbb{C}$, $m \in \mathbb{N}$ is a sequence of generalized paths such that $\diam(\gamma_m(X_m)) \to \infty$ as $m \to \infty$, then $\length(\gamma_m) \to 1$ as $m \to \infty$.

On the other hand, if we define $\gamma_m: [0,1] \to \mathbb{C}$ by $\gamma_m(t) = e^{2\pi imt}$, then $\length(\gamma_m) \to 1$ as $m \to \infty$, even though $\diam(\gamma_m([0,1])) = 2$ for all $m$.

\section{Properties of the function $\length$}
\label{sec:properties}

Let $n \geq 2$ be fixed.  All results in this section will be stated for $\length = \length_n$, and proofs will be given for the case $n = 2$.

The following basic properties follow immediately from the definition of the function $\length$.

\begin{prop}
\label{prop:basic properties 1}
Let $\gamma: X \to \mathbb{R}^n$ be a generalized path.
\begin{enumerate}[label=(\roman{*})]
\item \label{enum:zero length} $\length(\gamma) = 0$ if and only if $\gamma$ is a constant function.
\item \label{enum:homeo inv} If $h: Y \to X$ is a homeomorphism, then $\length(\gamma \circ h) = \length(\gamma)$.
\item \label{enum:iso inv} If $\Phi: \mathbb{R}^n \to \mathbb{R}^n$ is an isometry, then $\length(\Phi \circ \gamma) = \length(\gamma)$.
\end{enumerate}
\end{prop}

For the next properties, we need to consider a more restricted class of locally connected continua, namely dendrites.  A \emph{dendrite} is a locally connected continuum which contains no simple closed curve.  A characteristic feature of dendrites is that they are hereditarily unicoherent; that is, given any two intersecting subcontinua $A$ and $B$ of a dendrite $X$, the intersection $A \cap B$ is connected.  See Section \ref{sec:generalized paths} for examples to illustrate how these properties can fail when the domain of a generalized path is not a dendrite.

\begin{prop}
\label{prop:basic properties 2}
Let $X$ be a dendrite, and let $\gamma: X \to \mathbb{R}^n$ be a generalized path.
\begin{enumerate}[label=(\roman{*})]
\item \label{enum:monotone} If $A$ is a subcontinuum of $X$, then $\length(\gamma {\upharpoonright}_A) \leq \length(\gamma)$.  Moreover, $\length(\gamma {\upharpoonright}_A) = \length(\gamma)$ if and only if $\gamma$ is constant on each component of $X \smallsetminus A$.
\item \label{enum:subadditive} If $A,B$ are subcontinua of $X$ with $A \cup B = X$, then
\[ \length(\gamma) \leq \length(\gamma {\upharpoonright}_A) + \length(\gamma {\upharpoonright}_B) .\]
\end{enumerate}
\end{prop}

\begin{proof}
We treat the case $n = 2$.

Fix $x,t,\mu$, and for convenience denote $S^{x,t,\mu}_j$ and $C^{x,t,\mu}_n$ (defined as in Section \ref{sec:definition}) simply by $S_j$ and $C_n$, respectively.

Let $A \subseteq X$ be a subcontinuum.  Given $j \in \mathbb{Z}$ and a component $C$ of $(\gamma {\upharpoonright}_A)^{-1}(S_j)$, there exists some $n$ such that $C \subseteq C_n$.  Since $C_n \cap A$ is connected (by hereditary unicoherence), it follows that $C = C_n \cap A$.

Therefore there exists an injective function $\sigma$ from the non-negative integers to themselves such that ${\langle C_{\sigma(n)} \cap A \rangle}_{n=0}^\infty$ enumerates the collection of all components of the sets $(\gamma {\upharpoonright}_A)^{-1}(S^{x,t,\mu}_j)$ ($j \in \mathbb{Z}$) which have non-degenerate image under the map $\proj^\perp_t$, so that $\|\gamma(C_{\sigma(n)} \cap A)\|_t \geq \|\gamma(C_{\sigma(n+1)} \cap A)\|_t$ for all $n$.  Then
\begin{align*}
L^{x,t,\mu}(\gamma {\upharpoonright}_A) &= \sum_{n=0}^\infty \frac{\|\gamma(C_{\sigma(n)} \cap A)\|_t}{2^n} \\
&\leq \sum_{n=0}^\infty \frac{\|\gamma(C_{\sigma(n)})\|_t}{2^n} \\
&\leq L^{x,t,\mu}(\gamma) \textrm{\qquad (by the observation ($\ast$)).}
\end{align*}

Since this holds for all $x,t,\mu$, we have established the first statement of \ref{enum:monotone}.

For the second statement of \ref{enum:monotone}, suppose $\gamma$ is non-constant on some component $K$ of $X \smallsetminus A$.  The intersection $\overline{K} \cap A$ consists of a single point (see e.g.\ 10.9 and 10.24 of \cite{Nadler1992}).  Let $\{p\} = \overline{K} \cap A$, and let $q \in K$ be such that $\gamma(p) \neq \gamma(q)$.  There is a positive measure set of parameters $x,t,\mu$ and an integer $j \in \mathbb{Z}$ for which $\gamma(q) \in \interior(S^{x,t,\mu}_j)$ and $\gamma(p) \notin S^{x,t,\mu}_j$.  For such $x,t,\mu,j$, there is a component of $\gamma^{-1}(S^{x,t,\mu}_j)$ contained in $K$, which contributes positively to the sum $L^{x,t,\mu}(\gamma)$, thereby making it larger than $L^{x,t,\mu}(\gamma {\upharpoonright}_A)$.  It follows that $\length(\gamma) > \length(\gamma {\upharpoonright}_A)$.  The converse implication is immediate.

\medskip
Now suppose $A,B \subseteq X$ are subcontinua with $A \cup B = X$.  As above, for any $j \in \mathbb{Z}$, each component of $(\gamma {\upharpoonright}_A)^{-1}(S_j)$ (respectively $(\gamma {\upharpoonright}_B)^{-1}(S_j)$) has the form $C_n \cap A$ (respectively $C_n \cap B$) for some $n$.

Let ${\langle n(\alpha) \rangle}_{\alpha=0}^\infty$ and ${\langle m(\beta) \rangle}_{\beta=0}^\infty$ be the strictly increasing sequences of non-negative integers such that ${\langle C_{n(\alpha)} \cap A \rangle}_{\alpha=0}^\infty$ enumerates the collection of all components of the sets $(\gamma {\upharpoonright}_A)^{-1}(S^{x,t,\mu}_j)$ ($j \in \mathbb{Z}$) which have non-degenerate image under the map $\proj^\perp_t \circ \gamma$, and ${\langle C_{m(\beta)} \cap B \rangle}_{\beta=0}^\infty$ enumerates the collection of all components of the sets $(\gamma {\upharpoonright}_B)^{-1}(S^{x,t,\mu}_j)$ ($j \in \mathbb{Z}$) which have non-degenerate image under the map $\proj^\perp_t \circ \gamma$.  Note that these enumerations are not necessarily ordered according to the sizes of the images under $\proj^\perp_t \circ \gamma$.

For any $n$, we clearly have $\|\gamma(C_n)\|_t \leq \|\gamma(C_n \cap A)\|_t + \|\gamma(C_n \cap B)\|_t$.  Therefore
\begin{align*}
L^{x,t,\mu}(\gamma) &= \sum_{n=0}^\infty \frac{\|\gamma(C_n)\|_t}{2^n} \\
&\leq \sum_{n=0}^\infty \frac{\|\gamma(C_n \cap A)\|_t}{2^n} + \sum_{n=0}^\infty  \frac{\|\gamma(C_n \cap B)\|_t}{2^n} \\
&= \sum_{\alpha=0}^\infty \frac{\|\gamma(C_{n(\alpha)} \cap A)\|_t}{2^{n(\alpha)}} + \sum_{\beta=0}^\infty \frac{\|\gamma(C_{m(\beta)} \cap B)\|_t}{2^{m(\beta)}} \\
&\leq \sum_{\alpha=0}^\infty \frac{\|\gamma(C_{n(\alpha)} \cap A)\|_t}{2^\alpha} + \sum_{\beta=0}^\infty \frac{\|\gamma(C_{m(\beta)} \cap B)\|_t}{2^\beta} \textrm{\qquad (since $\alpha \leq n(\alpha)$, $\beta \leq m(\beta)$)} \\
&\leq L^{x,t,\mu}(\gamma {\upharpoonright}_{A}) + L^{x,t,\mu}(\gamma {\upharpoonright}_{B})  \textrm{\qquad (by the observation ($\ast$)).}
\end{align*}

Since this holds for all $x,t,\mu$, we have established \ref{enum:subadditive}.
\end{proof}

\begin{prop}
\label{prop:basic properties 3}
Let $z_1,z_2 \in \mathbb{R}^n$.  If $\gamma: X \to \mathbb{R}^n$ is any generalized path such that $z_1,z_2 \in \gamma(X)$, then $\length(\overline{z_1 z_2}) \leq \length(\gamma)$.  Moreover, if $\gamma(X)$ is not the straight line segment joining $z_1$ and $z_2$, or if $\gamma^{-1}(w)$ is disconnected for some $w$ on the straight line segment between $z_1$ and $z_2$, then $\length(\overline{z_1 z_2}) < \length(\gamma)$.
\end{prop}

Proposition \ref{prop:basic properties 3} can be proved directly from the definition of the function $\length$, and we leave this to the reader. Note that it also follows that if in a path $\gamma: [0,1] \to \mathbb{R}^n$ we replace the subpath $\gamma {\upharpoonright}_{[a,b]}$ with the straight line segment $\overline{\gamma(a) \gamma(b)}$ and if we denote the resulting path by $\gamma^*$, then $\length(\gamma^*) \leq \length(\gamma)$ with strict inequality if $\gamma {\upharpoonright}_{[a,b]}$ is not a monotone parametrization of the straight line segment $\overline{\gamma(a) \gamma(b)}$.

Next we consider $C(X) = C(X,\mathbb{R}^n)$, the set of all generalized paths $X \to \mathbb{R}^n$.  This is a metric space with the usual metric $\dsup(\gamma_1,\gamma_2) = \sup_{p \in X} |\gamma_1(p) - \gamma_2(p)|$.

\begin{prop}
\label{prop:length continuous}
The function $\length: C(X) \to \mathbb{R}^n$ is continuous.
\end{prop}

\begin{proof}
We treat the case $n = 2$.  Let $\gamma_0$ be in $C(X)$.

Suppose $\alpha < \length(\gamma_0) < \beta$.  We will prove that for small enough $\xi > 0$, if $\gamma \in C(X)$ with $\dsup(\gamma,\gamma_0) < \xi$, then $\alpha < \length(\gamma_0) < \beta$.

A simple modification of the proof of Lemma \ref{lem:lower semicontinuity} shows that for $\xi > 0$ small enough, if $\dsup(\gamma,\gamma_0) < \xi$ then $\length(\gamma) > \alpha$.  Thus it remains to show $\length(\gamma) < \beta$ for sufficiently small $\xi > 0$.

Fix a countable dense set $\{q_k\}_{k=1}^\infty \subset X$.  Given $k \neq l$ and $j \in \mathbb{Z}$, let
\begin{align*}
B_{kl}^j =& \{ (x,t,\mu) \in [0,1] \times [0,1] \times (0,1]: \\
& \qquad \textrm{there is a continuum } C \subseteq \gamma_0^{-1}(S^{x,t,\mu}_j) \textrm{ with } q_k,q_l \in C \textrm{ and} \\
& \qquad \textrm{for every such } C \textrm{ we have } \gamma_0(C) \cap \partial S^{x,t,\mu}_j \neq \emptyset \}
\end{align*}
and let $B = \bigcup_{\substack{k \neq l \\ j \in \mathbb{Z}}} B_{kl}^j$.  It is easy to see that $([0,1] \times [0,1] \times (0,1]) \smallsetminus B_{kl}^j$ is open, and so $B$ is $F_\sigma$, hence measurable.

\begin{claim}
\label{claim:bad measure zero}
$B$ has measure zero.
\end{claim}

\begin{proof}[Proof of Claim \ref{claim:bad measure zero}]
\renewcommand{\qedsymbol}{\textsquare (Claim \ref{claim:bad measure zero})}
Fix $k \neq l$ and $j \in \mathbb{Z}$.  It will be convenient to change variables from $(x,t,\mu)$ to $(z,t,\mu)$ so that for any fixed rotation angle $t$ and translation parameter $z$, as the strip width $\mu$ shrinks, the $j$-th strip itself shrinks inwards, nesting down on a line.

Given $(x,t,\mu) \in [0,1] \times [0,1] \times (0,1]$, let $z = \mu(x + \frac{1}{2} + j) \in (-\infty,\infty)$, and define $\Phi(x,t,\mu) = (z,t,\mu)$.

Observe that for $(z,t,\mu)$ in the image of $\Phi$, $\Phi^{-1}(z,t,\mu) = (\frac{z}{\mu} - \frac{1}{2} - j, t, \mu)$.  Thus $\Phi(B_{kl}^j)$ is contained in the set
\begin{align*}
B' =& \{(z,t,\mu): \textrm{there is a continuum } C \subseteq \gamma_0^{-1}(T^{z,t,\mu}) \textrm{ with } q_k,q_l \in C \textrm{ and} \\
& \qquad \textrm{for every such } C \textrm{ we have } \gamma_0(C) \cap \partial T^{z,t,\mu} \neq \emptyset \}
\end{align*}
where $T^{z,t,\mu} = \mu e^{t \pi i}(S_j + i(\frac{z}{\mu} - \frac{1}{2} - j)) = e^{t \pi i}(\mu(S_j - \frac{1}{2}i - ji) + iz)$.  Observe that the strip $T^{z,t,\mu}$ is centered about the line $e^{t \pi i}(\mathbb{R} + iz)$, and if $\mu' < \mu$, then $T^{z,t,\mu'}$ is contained in the interior of $T^{z,t,\mu}$.  Thus for any fixed $z,t$, there can be at most one $\mu$ for which $(z,t,\mu) \in B'$.  By Fubini's theorem, this implies $B'$ has measure zero.  Since $\Phi(B_{kl}^j) \subseteq B'$, we have that $\Phi(B_{kl}^j)$ has measure zero as well.

A straightforward calculation shows that $\Phi$ is a $C^1$-diffeomorphism on $[0,1] \times [0,1] \times (0,1]$ with Jacobian equal to $\mu$.  Thus by the change of variables theorem \cite[Theorem 2.47]{Folland1999}, the measure of $\Phi(B_{kl}^j)$ is equal to
\[ \iiint\limits_{B_{kl}^j} \! \mu \, dx \, dt \, d\mu .\]
Since $\mu > 0$ and $\Phi(B_{kl}^j)$ has measure zero, it follows that $B_{kl}^j$ has measure zero as well.  Since $B = \bigcup_{\substack{k \neq l \\ j \in \mathbb{Z}}} B_{kl}^j$, the Claim follows.
\end{proof}

\begin{claim}
\label{claim:continuity good set}
Given $(x_0,t_0,\mu_0) \in ([0,1] \times [0,1] \times (0,1]) \smallsetminus B$ and $\varepsilon > 0$, there exists $\delta > 0$ and $\xi_0 > 0$ such that if $|x - x_0| < \delta$, $|t - t_0| < \delta$, $|\mu - \mu_0| < \delta$, and $\dsup(\gamma,\gamma_0) < \xi_0$, then $L^{x,t,\mu}(\gamma) < L^{x_0,t_0,\mu_0}(\gamma_0) + \varepsilon$.
\end{claim}

\begin{proof}[Proof of Claim \ref{claim:continuity good set}]
\renewcommand{\qedsymbol}{\textsquare (Claim \ref{claim:continuity good set})}
For $j \in \mathbb{Z}$, let $S'_j$ denote the narrower (closed) strip obtained from $S^{x_0,t_0,\mu_0}_j$ by moving the boundary lines in towards the middle a distance of $\frac{\varepsilon}{20}$ each.

Let ${\langle C_n \rangle}_{n=0}^\infty$ enumerate the collection of all components of the sets $\gamma_0^{-1}(S^{x_0,t_0,\mu_0}_j)$ ($j \in \mathbb{Z}$) which have non-degenerate image under the map $\proj^\perp_t$, ordered so that $\|\gamma_0(C_n)\|_{t_0} \geq \|\gamma_0(C_{n+1})\|_{t_0}$ for all $n$.  For each $n$, let $j(n)$ be the integer such that $\gamma_0(C_n) \subset S^{x_0,t_0,\mu_0}_{j(n)}$.  By Lemma \ref{lem:components large projections}, there are only finitely many components $C_0,\ldots,C_N$ such that $\gamma_0(C_n)$ meets the narrower strip $S'_{j(n)}$, for $0 \leq n \leq N$.

Fix some $n$ with $0 \leq n \leq N$.  Let $U_1,\ldots,U_r$ be a finite cover of $C_n \cap \gamma_0^{-1}(S'_{j(n)})$ by connected open subsets of $X$ whose closures are mapped by $\gamma_0$ into the interior of $S^{x_0,t_0,\mu_0}_{j(n)}$.  Let $k$ be such that $q_k \in U_1$, and for each $2 \leq i \leq r$ let $l(i)$ be such that $q_{l(i)} \in U_i$.  Then for each $2 \leq i \leq r$, since $(x_0,t_0,\mu_0) \notin B_{k \, l(i)}^{j(n)}$ and $C_n$ is a continuum in $\gamma_0^{-1}(S_{j(n)}^{x_0,t_0,\mu_0})$ containing $q_k$ and $q_{l(i)}$, there exists a continuum $K_i$ containing $q_k$ and $q_{l(i)}$ which is mapped by $\gamma_0$ into the interior of the strip $S_{j(n)}^{x_0,t_0,\mu_0}$.  Let $C_n' = \overline{U_1} \cup \bigcup_{2 \leq i \leq r} (\overline{U_i} \cup K_i)$.  Then $C_n'$ is a continuum which is mapped by $\gamma_0$ into the interior of the strip $S_{j(n)}^{x_0,t_0,\mu_0}$ and such that $C_n \cap \gamma_0^{-1}(S'_j) \subseteq C_n' \subset C_n$.

Having done this for each $0 \leq n \leq N$, let $\delta > 0$ be small enough and let $\xi_0 > 0$ be small enough so that if $|x - x_0| < \delta$, $|t - t_0| < \delta$, $|\mu - \mu_0| < \delta$, and $\dsup(\gamma,\gamma_0) < \xi_0$, then for each $0 \leq n \leq N$ we have:

\begin{enumerate}[label=(\roman{*})]
\item \label{enum:int strip} $\gamma(C_n')$ is contained in the interior of the strip $S^{x,t,\mu}_{j(n)}$,
\item \label{enum:bound} $\|\gamma(C_n')\|_t < \|\gamma_0(C_n')\|_{t_0} + \frac{\varepsilon}{4}$, and
\item \label{enum:small between} if $A \subset X$ with $\gamma_0(A)$ contained in between two consecutive narrowed strips $S'_j$ and $S'_{j+1}$, then $\|\gamma(A)\|_t < \frac{\varepsilon}{8}$.
\end{enumerate}

Note that if $0 \leq n \leq N$ and if $C$ is the component of $\gamma^{-1}(S^{x,t,\mu}_{j(n)})$ containing $C_n'$, then $C$ consists of $C_n'$ plus some part which $\gamma_0$ maps in between $S'_{j(n)}$ and $S'_{j(n)-1}$, and some part which $\gamma_0$ maps in between $S'_{j(n)}$ and $S'_{j(n)+1}$.  Therefore, by \ref{enum:bound} and \ref{enum:small between} we have
\[ \|\gamma(C)\|_t < \|\gamma_0(C_n')\|_{t_0} + \tfrac{\varepsilon}{4} + 2 \cdot \tfrac{\varepsilon}{8} = \|\gamma_0(C_n')\|_{t_0} + \tfrac{\varepsilon}{2} .\]
Every other component $\tilde{C}$ of $\gamma^{-1}(S^{x,t,\mu}_j)$ satisfies $\|\gamma(\tilde{C})\|_t < \frac{\varepsilon}{8}$ by \ref{enum:small between}.  It follows that
\begin{align*}
L^{x,t,\mu}(\gamma) &< \sum_{n=0}^N \frac{\|\gamma_0(C_n')\|_{t_0} + \frac{\varepsilon}{2}}{2^n} + \sum_{n=N+1}^\infty \frac{\varepsilon/8}{2^n} \\
&< \sum_{n=0}^N \frac{\|\gamma_0(C_n')\|_{t_0}}{2^n} + \varepsilon \\
&\leq L^{x_0,t_0,\mu_0}(\gamma_0) + \varepsilon .
\end{align*}
\end{proof}

\bigskip
We are now ready to show that $\length(\gamma) < \beta$ for $\gamma$ sufficiently close to $\gamma_0$.

Recalling that $L^{x,t,\mu}(\gamma_0) < 2\mu \leq 2$, choose a step function $\psi = 2 - \sum_{i=0}^k c_i \,{\raisebox{1.5pt}{$\chi$}}_{A_i}$, where the $A_i \subset [0,1] \times [0,1] \times (0,1]$ are pairwise disjoint compact sets and ${\raisebox{1.5pt}{$\chi$}}_{A_i}$ is the characteristic function of the set $A_i$, with
\[ L^{x,t,\mu}(\gamma_0) \leq \psi(x,t,\mu)  \quad \textrm{for all } x,t,\mu \]
and
\[ \int_0^1 \! \int_0^1 \! \int_0^1 \! \psi(x,t,\mu) \, dx \, dt \, d\mu < \beta .\]
Let $\eta = \beta - \int_0^1 \! \int_0^1 \! \int_0^1 \! \psi \, dx \, dt \, d\mu > 0$.  By Claim \ref{claim:bad measure zero}, we can find a compact set $\Omega \subset ([0,1] \times [0,1] \times (0,1]) \smallsetminus B$ of measure $\geq 1 - \frac{\eta}{4}$.

Using Claim \ref{claim:continuity good set} and compactness of the sets $A_i \cap \Omega$, we can find $\xi_i$ small enough so that if $\dsup(\gamma,\gamma_0) < \xi_i$, then $L^{x,t,\mu}(\gamma) < \psi(x,t,\mu) + \frac{\eta}{4}$ for all $(x,t,\mu) \in A_i \cap \Omega$.  Letting $\xi = \min_i \xi_i$, it follows that if $\dsup(\gamma,\gamma_0) < \xi$, then
\begin{align*}
\length(\gamma) &= \int_0^1 \! \int_0^1 \! \int_0^1 \! L^{x,t,\mu}(\gamma) \, dx \, dt \, d\mu \\
&\leq \iiint\limits_\Omega \! L^{x,t,\mu}(\gamma) \, dx \, dt \, d\mu + 2 \cdot \tfrac{\eta}{4} \\
&\leq \iiint\limits_\Omega \! \left( \psi(x,t,\mu) + \tfrac{\eta}{4} \right) \, dx \, dt \, d\mu + 2 \cdot \tfrac{\eta}{4} \\
&\leq (\beta - \eta) + \tfrac{\eta}{4} + 2 \cdot \tfrac{\eta}{4} \\
&< \beta .
\end{align*}
\end{proof}

It follows from Proposition \ref{prop:length continuous} that for any path $\gamma_0$ and any $\varepsilon > 0$, there exists $\delta > 0$ such that if $\dsup(\gamma,\gamma_0) < \delta$, then $|\length(\gamma {\upharpoonright}_{[0,t]}) - \length(\gamma_0 {\upharpoonright}_{[0,t]})| < \varepsilon$ for all $t \in [0,1]$.

To see this, note that by Proposition \ref{prop:length continuous}, for any $t_0 \in [0,1]$, there is a small open interval $J_0$ around $t_0$ and $\delta_0 > 0$ small enough such that if $\dsup(\gamma,\gamma_0) < \delta_0$, then $|\length(\gamma {\upharpoonright}_{[0,t]}) - \length(\gamma_0 {\upharpoonright}_{[0,t_0]})| < \frac{\varepsilon}{2}$ for any $t \in J_0$.  Take a finite cover of $[0,1]$ by such intervals $J_0$ and take $\delta$ to be the minimum of the corresponding numbers $\delta_0$.  Suppose $\dsup(\gamma,\gamma_0) < \delta$.  Given any $t \in [0,1]$, take one of the intervals $J_0$ from the cover such that $t \in J_0$.  Then we have
\begin{align*}
|\length(\gamma {\upharpoonright}_{[0,t]}) - \length(\gamma_0 {\upharpoonright}_{[0,t]})| &\leq |\length(\gamma {\upharpoonright}_{[0,t]}) - \length(\gamma_0 {\upharpoonright}_{[0,t_0]})| + |\length(\gamma_0 {\upharpoonright}_{[0,t]}) - \length(\gamma_0 {\upharpoonright}_{[0,t_0]})| \\
&< \frac{\varepsilon}{2} + \frac{\varepsilon}{2} = \varepsilon .
\end{align*}

\bigskip
A consequence of Proposition \ref{prop:length continuous} is that for a path $\gamma: [a,b] \to \mathbb{R}^n$, $\length(\gamma)$ is small if and only if $\diam(\gamma([a,b]))$ is small.  This will suffice for our purposes, but in fact one can argue from the definition of $\length = \length_n$ that there are constants $c_1(n),c_2(n) > 0$ such that:
\begin{statement}{($\ast\ast$)}
If $\gamma: [a,b] \rightarrow \mathbb{R}^n$ is a path with $\diam(\gamma([a,b])) \leq \frac{1}{2}$, then
\[ c_1(n) \cdot \diam(\gamma([a,b])) \leq \length_n(\gamma) \leq c_2(n) \cdot \diam(\gamma([a,b])) .\]
\end{statement}

\section{Parameterization by $\length$}
\label{sec:parameterization length}

Let $n \geq 2$ be fixed.  As before, all results in this section will be stated for $\length = \length_n$, and proofs will be given for the case $n = 2$.

In this section, we work with $\mathcal{C}[0,1] = \mathcal{C}([0,1],\mathbb{R}^n)$, the set of all paths $\gamma: [0,1] \to \mathbb{R}^n$.  This is a metric space with the usual metric $\dsup(\gamma_1,\gamma_2) = \sup_{t \in [0,1]} |\gamma_1(t) - \gamma_2(t)|$.

\begin{defn}
\label{defn:reparam}
Given two paths $\gamma_1,\gamma_2: [0,1] \to \mathbb{R}^n$, we say that $\gamma_2$ is a \emph{reparameterization} of $\gamma_1$ if there are non-decreasing onto maps $m_1,m_2: [0,1] \to [0,1]$ such that $\gamma_i$ is constant on each fiber $m_i^{-1}(s)$, $s \in [0,1]$, for both $i = 1,2$, and $\gamma_1 \circ m_1^{-1} = \gamma_2 \circ m_2^{-1}$.  In this case, we write $\gamma_1 \approx \gamma_2$.
\end{defn}

Thus $\gamma_1 \approx \gamma_2$ if they both parameterize the same path, with the same orientation, where we disregard any constant sections.  Note that if $\gamma_1 \approx \gamma_2$, then $\length(\gamma_1) = \length(\gamma_2)$.  It is easy to see that $\approx$ is an equivalence relation on $\mathcal{C}[0,1]$.  Denote by $[\gamma]$ the equivalence class of $\gamma$ with respect to $\approx$.

Let $\Pi$ denote the collection of all equivalence classes $[\gamma]$.  We define a metric $\rho$ on $\Pi$ as follows:
\[ \rho([\gamma_1],[\gamma_2]) = \inf \{\sup_{t \in [0,1]} |\lambda_1(t) - \lambda_2(t)|: \lambda_1 \in [\gamma_1], \lambda_2 \in [\gamma_2]\} .\]
In fact, by reparameterizing, this can be expressed as $\rho([\gamma_1],[\gamma_2]) = \inf \{\sup_{t \in [0,1]} |\lambda_1(t) - \gamma_2(t)|: \lambda_1 \in [\gamma_1]\}$.  It is easy to show that $\rho$ is a metric, and that the resultant metric topology on $\Pi$ coincides with the quotient topology induced from $\mathcal{C}[0,1]$.

\bigskip
One can deduce from Propositions \ref{prop:basic properties 2}\ref{enum:monotone} and \ref{prop:length continuous} that given a path $\gamma: [0,1] \to \mathbb{R}^n$, the function $[0,1] \to [0,1)$ defined by $t \mapsto \length(\gamma {\upharpoonright}_{[0,t]})$ is continuous and non-decreasing.  As a result, we can make the following definition:

\begin{defn}
The \emph{standard parameterization} $\widetilde{\gamma}: [0,1] \to \mathbb{R}^n$ of $\gamma$, also called the parameterization of $\gamma$ by $\length$, is defined as follows.  If $\gamma$ is constant, then $\widetilde{\gamma} = \gamma$.  Otherwise, given $s \in [0,1]$, $\widetilde{\gamma}(s) = \gamma(t)$, where $t \in [0,1]$ is such that $\length(\gamma {\upharpoonright}_{[0,t]}) = s \cdot \length(\gamma)$.
\end{defn}
Note that this value $t$ may not be unique, but by Proposition \ref{prop:basic properties 2}\ref{enum:monotone}, the point $\gamma(t)$ is uniquely determined by $s$.  One can easily check that $\widetilde{\gamma}$ is a path (i.e.\ is a continuous function), $\widetilde{\gamma} \approx \gamma$, and $\length(\widetilde{\gamma} {\upharpoonright}_{[0,s]}) = s \cdot \length(\gamma)$ for any $s \in [0,1]$.  However, note that in general $\length(\widetilde{\gamma} {\upharpoonright}_{[s_1,s_2]}) \neq (s_2 - s_1) \length(\gamma)$ when $0 < s_1 < s_2 \leq 1$.

For the Euclidean path length, such a parameterization is only available for rectifiable paths, i.e.\ those paths with finite Euclidean length.

\bigskip
Observe that the standard parameterization is unique within each equivalence class of paths, in the sense that if $\gamma_1 \approx \gamma_2$, then $\widetilde{\gamma_1} = \widetilde{\gamma_2}$.

Consider the standard parameterization as a function $\Pi \to \mathcal{C}[0,1]$ which maps each class $[\gamma]$ to the unique standard parameterization $\widetilde{\gamma} \in [\gamma]$.  Denote by $\widetilde{\Pi}$ the range of this function; that is, $\widetilde{\Pi}$ is the set of all standard parameterizations of paths $[0,1] \to \mathbb{R}^n$.

\begin{thm}
\label{thm:Pi eq tildePi}
$\widetilde{\Pi}$ is a closed subset of $\mathcal{C}[0,1]$, and the function $[\gamma] \mapsto \widetilde{\gamma}$ is a homeomorphism from $\Pi$ to $\widetilde{\Pi}$.
\end{thm}

\begin{proof}
Suppose $\gamma \in \mathcal{C}[0,1] \smallsetminus \widetilde{\Pi}$, which means that $\length(\gamma {\upharpoonright}_{[0,s]}) \neq s \cdot \length(\gamma)$ for some $s \in [0,1]$.  Then for all $\lambda \in \mathcal{C}[0,1]$ which are uniformly close to $\gamma$, we have that $\lambda {\upharpoonright}_{[0,s]}$ is uniformly close to $\gamma {\upharpoonright}_{[0,s]}$ as well, hence by Proposition \ref{prop:length continuous} we have that $\length(\lambda)$ and $\length(\lambda {\upharpoonright}_{[0,s]})$ are close to $\length(\gamma)$ and $\length(\gamma {\upharpoonright}_{[0,s]})$, respectively.  It follows that $\length(\lambda {\upharpoonright}_{[0,s]}) \neq s \cdot \length(\lambda)$ if $\lambda$ is sufficiently close to $\gamma$, hence $\lambda \notin \widetilde{\Pi}$.  Thus $\mathcal{C}[0,1] \smallsetminus \widetilde{\Pi}$ is open, and so $\widetilde{\Pi}$ is closed.

It is clear that $[\gamma] \mapsto \widetilde{\gamma}$ is one-to-one, and that the inverse of this map is continuous, by definition of the metric $\rho$ on $\Pi$ (indeed the map $\widetilde{\gamma} \mapsto [\widetilde{\gamma}]$ is Lipschitz continuous with constant $1$).

To see that $[\gamma] \mapsto \widetilde{\gamma}$ is continuous, suppose $[\gamma_i]$ is a sequence in $\Pi$ converging to $[\gamma_\infty] \in \Pi$ (in the metric $\rho$ on $\Pi$).  By changing representatives if necessary, we may assume that $\gamma_i \to \gamma_\infty$ uniformly.  By Proposition \ref{prop:length continuous} (and the statements immediately after), it follows that for every $\varepsilon > 0$ there exists $n_0$ such that for all $i \geq n_0$ and all $t \in [0,1]$, $|\length(\gamma_i {\upharpoonright}_{[0,t]}) - \length(\gamma_\infty {\upharpoonright}_{[0,t]})| < \varepsilon$.

Fix $\varepsilon > 0$.  Let $\delta > 0$ be small enough so that for all $i \geq 1$ and all $t_1,t_2 \in [0,1]$, if $|\length(\gamma_i {\upharpoonright}_{[0,t_1]}) - \length(\gamma_i {\upharpoonright}_{[0,t_2]})| < \delta$ then $\diam(\gamma_i([t_1,t_2])) < \frac{\varepsilon}{2}$.  Let $n_0$ be large enough so that for all $i \geq n_0$ and $t \in [0,t]$, $|\length(\gamma_i {\upharpoonright}_{[0,t]}) - \length(\gamma_\infty {\upharpoonright}_{[0,t]})| < \frac{\delta}{2}$ and $|\gamma_i(t) - \gamma_\infty(t)| < \frac{\varepsilon}{2}$.

Given $s \in [0,1]$ and $i \geq n_0$, let $t_i,t_\infty \in [0,1]$ be such that $\length(\gamma_i {\upharpoonright}_{[0,t_i]}) = s \cdot \length(\gamma_i)$ and $\length(\gamma_\infty {\upharpoonright}_{[0,t_\infty]}) = s \cdot \length(\gamma_\infty)$, so that $\widetilde{\gamma}_i(s) = \gamma_i(t_i)$ and $\widetilde{\gamma}_\infty(s) = \gamma_\infty(t_\infty)$.  We have
\begin{align*}
|\length(\gamma_i {\upharpoonright}_{[0,t_i]}) - \length(\gamma_i {\upharpoonright}_{[0,t_\infty]})| &\leq |\length(\gamma_i {\upharpoonright}_{[0,t_i]}) - \length(\gamma_\infty {\upharpoonright}_{[0,t_\infty]})| + |\length(\gamma_\infty {\upharpoonright}_{[0,t_\infty]}) - \length(\gamma_i {\upharpoonright}_{[0,t_\infty]})| \\
&= |s \cdot \length(\gamma_i) - s \cdot \length(\gamma_\infty)| + |\length(\gamma_\infty {\upharpoonright}_{[0,t_\infty]}) - \length(\gamma_i {\upharpoonright}_{[0,t_\infty]})| \\
&\leq s \cdot \frac{\delta}{2} + \frac{\delta}{2} \\
&\leq \delta .
\end{align*}
By the definition of $\delta$, it follows that $\diam(\gamma_i([t_i,t_\infty])) < \frac{\varepsilon}{2}$.  This implies
\begin{align*}
|\widetilde{\gamma}_i(s) - \widetilde{\gamma}_\infty(s)| &= |\gamma_i(t_i) - \gamma_\infty(t_\infty)| \\
&\leq |\gamma_i(t_i) - \gamma_i(t_\infty)| + |\gamma_i(t_\infty) - \gamma_\infty(t_\infty)| \\
&\leq \frac{\varepsilon}{2} + \frac{\varepsilon}{2} = \varepsilon .
\end{align*}
Thus $\widetilde{\gamma}_i \to \widetilde{\gamma}_\infty$ uniformly.  Therefore, $[\gamma] \mapsto \widetilde{\gamma}$ is continuous.
\end{proof}

Given a family $\mathcal{F} \subseteq \Pi$, define $\widetilde{\mathcal{F}} = \{\widetilde{\gamma}: [\gamma] \in \mathcal{F}\}$.

\begin{cor}
\label{cor:compact transfer}
A set $\mathcal{F} \subseteq \Pi$ is closed (respectively, compact) if and only if $\widetilde{\mathcal{F}}$ is a closed (respectively, compact) subset of $\mathcal{C}[0,1]$.
\end{cor}

\bigskip
A classical result from metric geometry (see e.g.\ \cite{Busemann1955}) is that if $L > 0$ and ${\langle \gamma_m \rangle}_{m=1}^\infty$ is a sequence of paths in a bounded set, with Euclidean path lengths $\leq L$, and if $\widetilde{\gamma}_m: [0,1] \to \mathbb{R}^n$ is the parameterization of $\gamma_m$ by Euclidean path length (with domain linearly rescaled to $[0,1]$), then the sequence ${\langle \widetilde{\gamma}_m \rangle}_{m=1}^\infty$ has a subsequence which converges uniformly to a path of finite Euclidean length.  This reparameterization is necessary, as standard examples show (consider e.g.\ $\gamma_m: [0,1] \to [0,1]$ defined by $\gamma_m(s) = s^m$).

We will now prove a version of this result for the function $\length$, where the uniform bound on length assumption is replaced by a weaker restriction on the number of long sections of the paths.  Moreover, we prove that this condition is in fact a characterization of those families of paths which can be parameterized so as to be equicontinuous.  A similar result is proved in \cite{Silverman1969} using Morse's length function.

\begin{thm}
\label{thm:equicontinuous}
Let $\mathcal{F} \subseteq \Pi$.  Suppose that
\begin{statement}{($\dagger$)} for each $\varepsilon > 0$, there is a positive integer $N$ such that for every $[\gamma] \in \mathcal{F}$, there is no collection of more than $N$ pairwise disjoint subintervals of $[0,1]$ whose images under $\gamma$ have diameters $\geq \varepsilon$.
\end{statement}
Then the family $\widetilde{\mathcal{F}} = \{\widetilde{\gamma}: [\gamma] \in \mathcal{F}\}$ is equicontinuous.

Conversely, if an equicontinuous family can be formed by choosing parameterizations of all the paths in $\mathcal{F}$, then $\mathcal{F}$ satisfies the property \textup{($\dagger$)}.
\end{thm}

\begin{proof}
We treat the case $n = 2$.  As usual, we identify $\mathbb{R}^2$ with $\mathbb{C}$.

Fix $\varepsilon > 0$.  Let $N \geq 1$ be such that for every $\gamma$ with $[\gamma] \in \mathcal{F}$, there is no collection of more than $N$ pairwise disjoint subintervals of $[0,1]$ whose images under $\gamma$ have diameters $\geq \frac{\varepsilon}{16}$.  Let $\delta = \frac{\varepsilon^2}{2^{N+7} \cdot N}$.

Suppose for a contradiction that for some $[\gamma] \in \mathcal{F}$ there exist $0 \leq s_1 < s_2 \leq 1$ with $s_2 - s_1 < \delta$ and $\rho(\widetilde{\gamma}(s_1), \widetilde{\gamma}(s_2)) \geq \varepsilon$.  Note that
\begin{align*}
\length(\widetilde{\gamma} {\upharpoonright}_{[0,s_2]}) &= s_2 \cdot \length(\gamma) \\
&< s_1 \cdot \length(\gamma) + \delta \cdot \length(\gamma) \\
&< \length(\widetilde{\gamma} {\upharpoonright}_{[0,s_1]}) + \delta .
\end{align*}

Let $t_0 \in [0,1]$ be such that the line $\{re^{t_0\pi i}: r \in \mathbb{R}\}$ is orthogonal to the segment $\overline{\widetilde{\gamma}(s_1) \widetilde{\gamma}(s_2)}$.  Define $W \subset [0,1] \times [0,1] \times (0,1]$ by
\[ W = [0,1] \times [t_0 - \tfrac{1}{4}, t_0 + \tfrac{1}{4}] \times [\tfrac{\varepsilon}{8}, \tfrac{\varepsilon}{4}] ,\]
where the interval $[t_0 - \frac{1}{4}, t_0 + \frac{1}{4}]$ should be considered reduced $\mod\ 1$ (i.e.\ it represents the set of all $t \in [0,1]$ such that one of $|t-t_0|$, $|t-(t_0-1)|$, or $|t-(t_0+1)|$ is $\leq \frac{1}{4}$).  Note that for any $(x,t,\mu) \in W$, any strip $S^{x,t,\mu}_j$ ($j \in \mathbb{Z}$) covers less than half of the line segment $\overline{\widetilde{\gamma}(s_1) \widetilde{\gamma}(s_2)}$.

Consider a fixed $x,t,\mu \in W$.  Let $C_0,\ldots,C_N$ and $D_0,\ldots,D_N$ be the first $N+1$ components of $(\widetilde{\gamma} {\upharpoonright}_{[0,s_1]})^{-1} (S^{x,t,\mu}_j)$ and $(\widetilde{\gamma} {\upharpoonright}_{[0,s_2]})^{-1} (S^{x,t,\mu}_j)$ ($j \in \mathbb{Z}$), respectively, ordered so that $\|\widetilde{\gamma}(C_i)\|_t \geq \|\widetilde{\gamma}(C_{i+1})\|_t$ and $\|\widetilde{\gamma}(D_i)\|_t \geq \|\widetilde{\gamma}(D_{i+1})\|_t$ for each $i = 0,1,\ldots,N-1$.  So $\sum_{i=0}^N \frac{\|\widetilde{\gamma}(C_i)\|_t}{2^i}$ and $\sum_{i=0}^N \frac{\|\widetilde{\gamma}(D_i)\|_t}{2^i}$ are the first $N+1$ terms of the sums $L^{x,t,\mu}(\widetilde{\gamma} {\upharpoonright}_{[0,s_1]})$ and $L^{x,t,\mu}(\widetilde{\gamma} {\upharpoonright}_{[0,s_2]})$, respectively.

Note that
\begin{equation}
\setcounter{equation}{1}
\|\widetilde{\gamma}(D_i)\|_t \geq \|\widetilde{\gamma}(C_i)\|_t \textrm{ for each } i = 0,1,\ldots,N .
\end{equation}
Moreover, there is some $j \in \{0,1,\ldots,N-1\}$ such that $D_j \subset (s_1,s_2)$ and $\|\widetilde{\gamma}(D_j)\|_t = \mu$.  Since such a component is absent in the list $C_0,\ldots,C_N$, we have $\|\widetilde{\gamma}(D_{i+1})\|_t \geq \|\widetilde{\gamma}(C_i)\|_t$ for each $i = j,\ldots,N-1$.

Now $\|\widetilde{\gamma}(D_j)\|_t = \mu \geq \frac{\varepsilon}{8}$, and $\|\widetilde{\gamma}(D_N)\|_t < \frac{\varepsilon}{16}$ by choice of $N$, so there must be some $i$ between $j$ and $N-1$ such that $\|\widetilde{\gamma}(D_i)\|_t > \|\widetilde{\gamma}(D_{i+1})\|_t + \frac{\varepsilon}{8N}$.  Hence
\begin{equation}
\setcounter{equation}{2}
\|\widetilde{\gamma}(D_i)\|_t > \|\widetilde{\gamma}(C_i)\|_t + \frac{\varepsilon}{8N} .
\end{equation}

It follows from (1) and (2) that
\begin{align*}
L^{x,t,\mu}(\widetilde{\gamma} {\upharpoonright}_{[0,s_2]}) &> L^{x,t,\mu} (\widetilde{\gamma} {\upharpoonright}_{[0,s_1]}) + \frac{\varepsilon/8N}{2^i} \\
&> L^{x,t,\mu}(\widetilde{\gamma} {\upharpoonright}_{[0,s_1]}) + \frac{\varepsilon/8N}{2^N} \\
&= L^{x,t,\mu}(\widetilde{\gamma} {\upharpoonright}_{[0,s_1]}) + \frac{\varepsilon}{2^{N+3} \cdot N}
\end{align*}

Noting that the measure of $W$ is $1 \cdot \frac{1}{2} \cdot (\frac{\varepsilon}{4} - \frac{\varepsilon}{8}) = \frac{\varepsilon}{16}$, it follows that
\begin{align*}
\length(\widetilde{\gamma} {\upharpoonright}_{[0,s_2]}) &\geq \length(\widetilde{\gamma} {\upharpoonright}_{[0,s_1]}) + \frac{\varepsilon}{2^{N+3} \cdot N} \cdot \frac{\varepsilon}{16} \\
&= \length(\widetilde{\gamma} {\upharpoonright}_{[0,s_1]}) + \delta .
\end{align*}
But this contradicts the assumption that $s_2 - s_1 < \delta$.

Thus for every $[\gamma] \in \mathcal{F}$, if $0 \leq s_1 < s_2 \leq 1$ with $s_2 - s_1 < \delta$, then $\rho(\widetilde{\gamma}(s_1), \widetilde{\gamma}(s_2)) < \varepsilon$.

\bigskip
For the converse, suppose there is some $\varepsilon > 0$ such that for any positive integer $N$, there exists a path $\gamma_N$ with $[\gamma_N] \in \mathcal{F}$ and a collection of $N$ disjoint subintervals of $[0,1]$ whose images under $\gamma_N$ have diameters $\geq \varepsilon$.  Note that at least one of these subintervals must have width $\leq \frac{1}{N}$; denote it by $J_N$.

Let $s \in [0,1]$ be an accumulation point of the centers of the intervals $J_N$, $N = 1,2,3,\ldots$.  Then for any $\delta > 0$, there is some $N$ such that $J_N \subset (s-\delta, s+\delta)$, and hence $\gamma_N((s-\delta, s+\delta))$ has diameter $\geq \varepsilon$.  Thus we cannot choose parameterizations of the paths in $\mathcal{F}$ to obtain an equicontinuous family.
\end{proof}

Theorem \ref{thm:equicontinuous} implies in particular that if it is possible to parameterize the paths of a family $\mathcal{F}$ to obtain an equicontinuous family, then the standard parameterization will accomplish this.

\begin{thm}
\label{thm:compact characterization}
Let $\mathcal{F} \subseteq \Pi$.  Then $\overline{\mathcal{F}}$ is compact if and only if the following two properties are satisfied:
\begin{enumerate}[label=(\arabic{*})]
\item \label{enum:initial bounded} the set $\{\gamma(0): [\gamma] \in \mathcal{F}\}$ is bounded; and
\item \label{enum:equicontinuous} $\mathcal{F}$ satisfies the property \textup{($\dagger$)} (from Theorem \ref{thm:equicontinuous}).
\end{enumerate}
\end{thm}

In particular, if $\mathcal{F}$ satisfies properties \ref{enum:initial bounded} and \ref{enum:equicontinuous}, then the closure of $\widetilde{\mathcal{F}}$ in $\mathcal{C}[0,1]$ is compact.

\begin{proof}
By the Arzel\`{a}-Ascoli theorem \cite[Theorem 4.43]{Folland1999}, the closure of $\widetilde{\mathcal{F}}$ is compact if and only if $\widetilde{\mathcal{F}}$ is equicontinuous and pointwise bounded, i.e.\ for every $t \in [0,1]$ the set $\{\widetilde{\gamma}(t): \widetilde{\gamma} \in \widetilde{\mathcal{F}}\}$ is bounded.  By Theorem \ref{thm:equicontinuous}, equicontinuity of $\widetilde{\mathcal{F}}$ is equivalent to $\mathcal{F}$ satisfying the property \textup{($\dagger$)}.  Moreover, in the presence of \textup{($\dagger$)}, the condition \ref{enum:initial bounded} is clearly equivalent to $\widetilde{\mathcal{F}}$ being pointwise bounded.

Finally, by Theorem \ref{thm:Pi eq tildePi}, $\overline{\mathcal{F}}$ is compact if and only if the closure of $\widetilde{\mathcal{F}}$ is compact.
\end{proof}


\section{Midpoint parameterization}
\label{sec:parameterization midpoint}

One drawback to the standard parameterization of a path, introduced in the previous section, is that it does not commute with reversal of orientation of a path.  That is, if we define $r: [0,1] \to [0,1]$ by $r(t) = 1-t$, then for a non-constant path $\gamma$, the standard parameterization of $\gamma \circ r$ is not the same as $\widetilde{\gamma} \circ r$.

In this section we introduce a second parameterization of a path and consider some applications, including a way to canonically define a homeomorphism between two arcs (or even between two families of arcs) once the endpoints have been assigned.

\bigskip
Given a path $\gamma$, we define the \emph{midpoint parameterization} $\gamma^\ast: [0,1] \to \mathbb{C}$ as follows.  Let $m \in (0,1)$ be such that $\length( \gamma {\upharpoonright}_{[0,m]} ) = \length( \gamma {\upharpoonright}_{[m,1]} ) = L$.  Define $\gamma_1, \gamma_2: [0,1] \to \mathbb{C}$ by $\gamma(t) = \gamma_1(m - mt)$ and $\gamma_2(t) = \gamma(m + (1-m)t)$, and consider their standard parameterizations $\widetilde{\gamma}_1$, $\widetilde{\gamma}_2: [0,1] \to \mathbb{C}$.  Then
\[ \gamma^\ast(t) =
\begin{cases}
\widetilde{\gamma}_1(1 - 2t) &\textrm{if } 0 \leq t \leq \frac{1}{2} \\
\widetilde{\gamma}_2(2t - 1) &\textrm{if } \frac{1}{2} < t \leq 1 \\
\end{cases}
\]
Observe that $\gamma^\ast \circ r = (\gamma \circ r)^\ast$.

As with the standard parameterization, the midpoint parameterization is unique within each equivalence class of paths, in the sense that if $\gamma_1 \approx \gamma_2$, then $\gamma_1^\ast = \gamma_2^\ast$.

We leave it to the reader to see that the following analogue of Theorem \ref{thm:Pi eq tildePi} also holds (where the standard parameterization is replaced by the midpoint paramterization).  We denote by $\Pi^\ast$ the set of all midpoint parameterizations of paths $[0,1] \to \mathbb{R}^n$.  Recall that the topology on $\Pi$ is given by the metric $\rho$ and that $\Pi^\ast$ is a subspace of $\mathcal{C}[0,1] = \mathcal{C}([0,1], \mathbb{R}^n)$ with the metric $\dsup(\gamma_1,\gamma_2) = \sup_{t \in [0,1]} |\gamma_1(t) - \gamma_2(t)|$.

\begin{thm}
\label{thm:Pi eq astPi}
$\Pi^\ast$ is a closed subset of $\mathcal{C}[0,1]$, and the function $[\gamma] \mapsto \gamma^\ast$ is a homeomorphism from $\Pi$ to $\Pi^\ast$.
\end{thm}

\bigskip
An \emph{arc} is a space $A$ which is homeomorphic to the interval $[0,1]$.  By a \emph{parametrization of an arc $A$} we mean a homeomorphism $\gamma: [0,1] \to A$.

Given two arcs $A_1$ and $A_2$, and a bijection $f$ between their endpoint sets, we can extend $f$ to a canonical homeomorphism $F: A_1 \to A_2$ as follows.  Choose a parameterization $\gamma_1$ of $A_1$, and a parameterization $\gamma_2$ of $A_2$ such that $\gamma_2(0) = f(\gamma_1(0))$ and $\gamma_2(1) = f(\gamma_1(1))$.  Then given $x \in A_1$, define $F(x) = \gamma_2^\ast((\gamma_1^\ast)^{-1}(x))$.  Observe that $F$ is independent of the choice of orientation of $\gamma_1$, since $\gamma_1^\ast \circ r = (\gamma_1 \circ r)^\ast$ and $\gamma_2^\ast \circ r = (\gamma_2 \circ r)^\ast$ (hence the word ``canonical'').

In the next result, we show that this canonical homeomorphism has some useful convergence properties.  

\begin{defn}
A collection of arcs $\lam_X$ is a \emph{lamination} of a space $X \subset \mathbb{R}^n$ if:
\begin{enumerate}[label=(\arabic{*})]
\item $X = \bigcup \lam_X$;
\item any two distinct arcs in $\lam_X$ meet at most in one common endpoint;
\item \label{enum:arc converge} given a sequence $\langle A_i \rangle_{i=1}^\infty$ of arcs in $\lam_X$:
\begin{enumerate}[label=(\alph{*})]
\item if $\diam A_i \to 0$, then any accumulation point of the arcs $A_i$ is either an endpoint of some arc in $\mathcal{L}_X$, or not in $X$;
\item otherwise, there is an arc $A_\infty \in \mathcal{L}_X$, a subsequence $\langle A_{i_j} \rangle_{j=1}^\infty$ and homeomorphisms $h_j: A_\infty \to A_{i_j}$ such that $\dsup(h_j,\id_{A_\infty}) \to 0$ as $j \to \infty$.
\end{enumerate}
\end{enumerate}

Denote by $\mathcal{E}(\lam_X)$ the set of all endpoints of the arcs in $\lam_X$.
\end{defn}

Note that the conclusion in condition \ref{enum:arc converge}(b) is equivalent to the statement that $A_\infty \cup \bigcup_j A_{i_j}$ is homeomorphic to the product of $[0,1]$ and the convergent sequence $\{0\} \cup \{\frac{1}{n}: n = 1,2,\ldots\}$.

\begin{thm}
\label{thm:lamination}
Let $\lam_X$ and $\lam_Y$ be laminations of spaces $X$ and $Y$, respectively.  Suppose $f: \mathcal{E}(\lam_X) \to \mathcal{E}(\lam_Y)$ is a continuous function which maps the endpoints of any arc in $\lam_X$ onto the set of endpoints of some arc in $\lam_Y$.  Then there exists a continuous extension $F: X \to Y$ of $f$ which is one-to-one on each arc in $\lam_X$.

Moreover, if additionally $f$ is a homeomorphism and $f^{-1}$ maps the endpoints of any arc in $\lam_Y$ to the endpoints of some arc in $\lam_X$, then $F$ is a homeomorphism.
\end{thm}

\begin{proof}
Let $\mathcal{P} = \{[\gamma]: \gamma$ parameterizes some arc $A \in \lam_X\}$ and $\mathcal{Q} = \{[\lambda]: \lambda$ parameterizes some arc $B \in \lam_Y\}$.  Define the function $\mathbf{g}: \mathcal{P} \to \mathcal{Q}$ by $\mathbf{g}([\gamma]) = [\lambda]$ if $f(\gamma(0)) = \lambda(0)$ and $f(\gamma(1)) = \lambda(1)$.

\begin{claim}
\label{claim:g continuous}
$\mathbf{g}$ is continuous.
\end{claim}

\begin{proof}[Proof of Claim \ref{claim:g continuous}]
\renewcommand{\qedsymbol}{\textsquare (Claim \ref{claim:g continuous})}
Suppose that $\langle [\gamma_i] \rangle_{i=1}^\infty$ is a sequence of elements of $\mathcal{P}$ converging to  $[\gamma_\infty] \in \mathcal{P}$.  By continuity of $f$, $\lim_{i \to \infty} f(\gamma_i(0)) = f(\gamma_\infty(0))$ and $\lim_{i \to \infty} f(\gamma_i(1)) = f(\gamma_\infty(1))$, and it follows from property \ref{enum:arc converge} (for $\lam_Y$) that $\lim_{i \to \infty} \mathbf{g}([\gamma_i]) = \mathbf{g}([\gamma_\infty])$.
\end{proof}

Define $F: X \to Y$ as follows.  Given $A \in \lam_X$, choose $\gamma$ parameterizing $A$.  Let $\gamma^\ast$ be the midpoint parameterization of $\gamma$, and let $\lambda^\ast$ be the midpoint parameterization of $\mathbf{g}([\gamma])$.  Now for each $x \in A$, define $F(x) = \lambda^\ast((\gamma^\ast)^{-1}(x))$.

Observe that the definition of $F$ on an arc $A$ does not depend on the choice of orientation of the parameterization $\gamma$ of $A$, because $\gamma^\ast \circ r = (\gamma \circ r)^\ast$, $\lambda^\ast \circ r = (\lambda \circ r)^\ast$, and $\mathbf{g}([\gamma \circ r]) = [\lambda \circ r]$.  Thus $F$ is well-defined on each arc $A$.  Moreover, if $x \in \mathcal{E}(\lam_X)$, then $F(x) = f(x)$.  Thus $F$ is well-defined on $X$ (since two arcs in $\lam_X$ meet at most in an endpoint), and extends $f$.  It is also clear that $F$ is one-to-one on any arc in $\lam_X$ since $\gamma^*$ and $\lambda^*$ are homeomorphisms.

\begin{claim}
\label{claim:F continuous}
$F$ is continuous.
\end{claim}

\begin{proof}[Proof of Claim \ref{claim:F continuous}]
\renewcommand{\qedsymbol}{\textsquare (Claim \ref{claim:F continuous})}
Suppose $\langle x_i \rangle_{i=1}^\infty$ is a sequence in $X$ converging to $x_\infty \in X$.  For each $n \in \mathbb{N} \cup \{\infty\}$, let $A_n \in \lam_X$ be an arc containing $x_n$, and let $B_n = F(A_n) \in \lam_Y$.

If $\diam(A_i) \to 0$, then by \ref{enum:arc converge} $x_\infty$ is an endpoint of $A_\infty$.  For each $i$, let $x_i'$ be an endpoint of $A_i$.  Then $x_i' \to x_\infty$, so by continuity of $f$, $f(x_i') \to f(x_\infty)$.  Thus both endpoints of the arcs $B_i$ converge to $f(x_\infty)$, so again by \ref{enum:arc converge} we have $\diam(B_i) \to 0$ and $f(x_i) \to f(x_\infty)$.

Otherwise, we may assume (by taking a subsequence of $x_i$, and rechoosing $A_\infty$ if necessary in the case that $x_\infty$ is an endpoint), that there are homeomorphisms $h_i: A_\infty \to A_i$ such that $\dsup(h_i,\id_{A_\infty}) \to 0$.  For each $n \in \mathbb{N} \cup \{\infty\}$, let $\gamma_n$ parameterize $A_n$, and let $\mathbf{g}([\gamma_n]) = [\lambda_n]$, so that $\lambda_n$ parameterizes $B_n$.  We may assume (by choosing appropriate orientations) that $\gamma_i(0) \to \gamma_\infty(0)$ and $\gamma_i(1) \to \gamma_\infty(1)$.  Then $[\gamma_i] \to [\gamma_\infty]$, and by continuity of $\mathbf{g}$, $[\lambda_i] \to [\lambda_\infty]$.

By Theorem \ref{thm:Pi eq astPi}, we have $\gamma_i^\ast \to \gamma_\infty^\ast$ and $\lambda_i^\ast \to \lambda_\infty^\ast$ uniformly.  For each $n \in \mathbb{N} \cup \{\infty\}$, let $t_n \in [0,1]$ be such that $\gamma_n(t_n) = x_n$.  Then by uniform convergence and continuity, $t_i \to t_\infty$, and $\lambda_i^\ast(t_i) \to \lambda_\infty^\ast(t_\infty)$.  Thus
\[ F(x_i) = \lambda_i^\ast((\gamma_i^\ast)^{-1}(x_i)) = \lambda_i^\ast(t_i) \to \lambda_\infty^\ast(t_\infty) = F(x_\infty) \]
as needed.
\end{proof}

If, in addition, $f$ is a homeomorphism and $f^{-1}$ also maps the endpoints of any arc in $\lam_Y$ to the endpoints of some arc in $\lam_X$, then $f^{-1}$ extends in the same way to a continuous function $Y \to X$ which is the inverse of $F$.  Thus $F$ is a homeomorphism.
\end{proof}

For a fixed space $X \subset \mathbb{R}^n$ with a lamination $\mathcal{L}_X$, let
\begin{align*}
M = \{(Y,\mathcal{L}_Y,f): & \textrm{ $Y \subset \mathbb{R}^n$ is bounded, $\mathcal{L}_Y$ is a lamination of $Y$, and} \\
& \textrm{$f:\mathcal{E}(\mathcal{L}_X) \to \mathcal{E}(\mathcal{L}_Y)$ satisfies the hypotheses of Theorem \ref{thm:lamination}}\} .
\end{align*}
Theorem \ref{thm:lamination} affords an operator $\Theta$ from $M$ to the set $\mathcal{C}_b(X,\mathbb{R}^n)$ of bounded continuous functions from $X$ into $\mathbb{R}^n$, where if $F = \Theta(Y,\mathcal{L}_Y,f)$ then $F(X) \subset Y$ and $F {\upharpoonright}_{\mathcal{E}(\mathcal{L}_X)} = f$.

We next prove that this operator is continuous, in the sense that if $\mathcal{L}_{Y_1}$ and $\mathcal{L}_{Y_2}$ are nearby laminations of two spaces $Y_1$ and $Y_2$, and if $f_1:\mathcal{E}(\mathcal{L}_X) \to \mathcal{E}(\mathcal{L}_{Y_1})$ and $f_2:\mathcal{E}(\mathcal{L}_X) \to \mathcal{E}(\mathcal{L}_{Y_2})$ are close functions as in Theorem \ref{thm:lamination}, then the extensions $F_1:X \to Y_1$ and $F_2:X \to Y_2$ are close as well.  To make this precise, we define a metric $\mathbf{d}$ on $M$ by
\begin{align*}
\mathbf{d} \left( (Y_1,\mathcal{L}_{Y_1},f_1), (Y_2,\mathcal{L}_{Y_2},f_2) \right) &= \sup \{\rho([\lambda_1],[\lambda_2]): \gamma \textrm{ parameterizes an arc $A \in \mathcal{L}_X$} \\
&\textrm{and } \lambda_i \textrm{ parameterizes the corresponding arc in $\mathcal{L}_{Y_i}$} \\
&\textrm{with $\lambda_i(0) = f_i(\gamma(0))$ and $\lambda_i(1) = f_i(\gamma(1))$ for $i=1,2$} \} .
\end{align*}
It is straightforward to verify that $\mathbf{d}$ is a metric on $M$.  On $\mathcal{C}_b(X,\mathbb{R}^n)$, we use the metric $\dsup$.

\begin{thm}
\label{thm:lamination cont}
Let $X \subset \mathbb{R}^n$, and let $\mathcal{L}_X$ be a lamination of $X$.  The operator $\Theta: M \to \mathcal{C}_b(X,\mathbb{R}^n)$ given by Theorem \ref{thm:lamination} is continuous.  Moreover, $\Theta(X, \mathcal{L}_X, \id_{\mathcal{E}(\mathcal{L}_X)}) = \id_X$.
\end{thm}

\begin{proof}
Let $(Y_0,\mathcal{L}_{Y_0},f_0) \in M$, and let $F_0 = \Theta(Y_0,\mathcal{L}_{Y_0},f_0) \in \mathcal{C}_b(X,\mathbb{R}^n)$.  Let $\varepsilon > 0$.

Since $Y_0$ is bounded, by the definition of a lamination, the set $\{[\lambda]: \lambda$ parameterizes an arc $B_0 \in \mathcal{L}_{Y_0}$ with $\diam(B_0) \geq \frac{\varepsilon}{2}\}$ is compact in $\Pi$.  Therefore, by Theorem \ref{thm:Pi eq astPi}, there exists $\delta > 0$ such that if $\lambda_0$ parameterizes an arc in $\mathcal{L}_{Y_0}$ of diameter $\geq \frac{\varepsilon}{2}$ and if $\lambda$ is any path with $\rho([\lambda_0],[\lambda]) < \delta$, then $\dsup(\lambda_0^*,\lambda^*) < \varepsilon$.  We may assume that $\delta \leq \frac{\varepsilon}{2}$.

Observe that if $\diam(\lambda([0,1])) < \frac{\varepsilon}{2}$ and $\rho([\lambda_0],[\lambda]) < \frac{\varepsilon}{2}$, then every point in the range of $\lambda_0$ is within $\varepsilon$ of each point in the range of $\lambda$, hence $\dsup(\lambda_0^*,\lambda^*) < \varepsilon$ as well.  Thus in fact for any arc $B_0 \in \mathcal{L}_Y$, if $\lambda_0$ parameterizes $B_0$ and if $\lambda$ is any path with $\rho([\lambda_0],[\lambda]) < \delta$, then $\dsup(\lambda_0^*,\lambda^*) < \varepsilon$.

Let $(Y,\mathcal{L}_Y,f) \in M$ with $\mathbf{d}((Y_0,\mathcal{L}_{Y_0},f_0),(Y,\mathcal{L}_Y,f)) < \delta$, and let $F = \Theta(Y,\mathcal{L}_Y,f)$.

Let $A \in \mathcal{L}_X$, let $\gamma$ parameterize $A$, and let $\lambda_0$ and $\lambda$ parameterize the corresponding arcs $B_0 \in \mathcal{L}_{Y_0}$ and $B \in \mathcal{L}_Y$ with $\lambda_0(0) = f_0(\gamma(0))$, $\lambda_0(1) = f_0(\gamma(1))$, $\lambda(0) = f(\gamma(0))$, and $\lambda(1) = f(\gamma(1))$.  By definition of $\mathbf{d}$, we have $\rho([\lambda_0],[\lambda]) < \delta$, hence by the choice of $\delta$, $\dsup(\lambda_0^*,\lambda^*) < \varepsilon$.  Moreover, by the definition of $F$ and $F_0$ from Theorem \ref{thm:lamination}, $\dsup(F_0 {\upharpoonright}_A, F {\upharpoonright}_A) = \dsup(\lambda_0^*,\lambda^*)$.  Thus since $A$ was arbitrary, we have $\dsup(F_0,F) < \varepsilon$.

The second part of this Theorem is clear from the definition of $\Theta$.
\end{proof}

\section{Generalized paths}
\label{sec:generalized paths}

To see that the assumption that $X$ is a dendrite in Proposition \ref{prop:basic properties 2} is necessary, consider the identity function $\id_\mathbb{D}$ on the unit disk $\mathbb{D} \subset \mathbb{C}$ with boundary circle $\mathbb{S}^1$.  It is not difficult to see that $\length(\id_\mathbb{D}) < \length(\id_{\mathbb{S}^1})$.

\begin{figure}
\begin{center}
\includegraphics{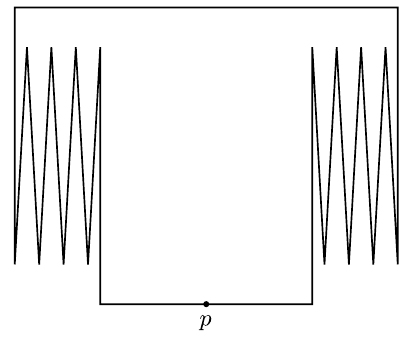}
\end{center}

\caption{A particular embedding of the circle in the plane.}
\label{fig:non dendrite}
\end{figure}

Moreover, consider the embedding $O$ of the circle $\mathbb{S}^1$ depicted in Figure \ref{fig:non dendrite}.  Let $\gamma: [0,1] \to O$ be a path which goes exactly once around the circle $O$, starting and ending at the indicated point $p$, and otherwise one-to-one.  We claim that $\length(\gamma) < \length(\id_O)$, which can be argued as follows:

Given a strip $S^{x,t,\mu}_j$ containing the point $p$, the component $C$ of $O \cap S^{x,t,\mu}_j$ containing $p$ corresponds to two components $[0,c]$ and $[d,1]$ of $\gamma^{-1}(S^{x,t,\mu}_j)$.  For nearly horizontal strips (i.e.\ for $t$ values close to $0$ or $1$) the sets $\proj^\perp_t(\gamma([0,c]))$ and $\proj^\perp_t(\gamma([d,1]))$ may overlap; however, because of the oscillation up and down on the left and right sides of the circle $O$, for such parameters $x,t,\mu$ there are many other components of $[0,1] \cap \gamma^{-1}(S^{x,t,\mu}_j)$ and of $O \cap S^{x,t,\mu}_j$ ($j \in \mathbb{Z}$) with large projections, hence the weighted sums $L^{x,t,\mu}(\id_O)$ and $L^{x,t,\mu}(\gamma)$ will differ only very slightly.  For all other values of $x,t,\mu$, the sets $\proj^\perp_t(\gamma([0,c]))$ and $\proj^\perp_t(\gamma([d,1]))$ share only the point $\proj^\perp_t(p)$, and one of them will be added with a smaller weight in the sum $L^{x,t,\mu}(\gamma)$ than that of $C$ in $L^{x,t,\mu}(\id_O)$.  In particular, this is so for values of $x,t,\mu$ for which the strips $S^{x,t,\mu}_j$ are wide and nearly vertical, and for these values resulting difference between $L^{x,t,\mu}(\id_O)$ and $L^{x,t,\mu}(\gamma)$ will be more pronounced due to the small number of terms in these sums.  Thus, with an appropriate amount of oscillation, we obtain that $\length(\gamma) < \length(\id_{O})$.

Now if we let $A$ be a very small arc in $O$ containing the point $p$ and such that $\length(\id_A) < \length(\id_O) - \length(\gamma)$, and let $A' = \overline{O \smallsetminus A}$, then it follows that $\length(\id_{O}) > \length(\id_A) + \length(\id_{A'})$.

\bibliographystyle{amsplain}
\bibliography{PathLength}

\end{document}